\documentclass[a4paper,12pt]{article}
\usepackage{amsfonts}
\usepackage{mathrsfs}
\usepackage{amssymb}
\usepackage{amsthm}
\usepackage{amsmath}
\usepackage{amscd}
\usepackage{graphicx}
\usepackage{xypic}
\usepackage[all,poly,knot]{xy}
\usepackage[dvipdfm,
            CJKbookmarks=true,
            bookmarksnumbered=true,
            bookmarksopen=true,
            colorlinks=true,
            citecolor=blue,
            linkcolor=blue,
            anchorcolor=green,
            urlcolor=blue
            ]{hyperref}

\usepackage{cite}

\newtheoremstyle{plain}{4pt}{4pt}{\itshape}{}{\bfseries}{}{2mm}{}
\theoremstyle{plain}
 \newtheorem{thm}{Theorem}[section]
 \newtheorem{lem}[thm]{Lemma\!}
 \newtheorem{exa}[thm]{Example\!}
 \newtheorem{defn}[thm]{Definition\!}
 \newtheorem{cor}[thm]{Corollary\!}
 \newtheorem{rem}[thm]{Remark\!}
 \newtheorem{pro}[thm]{Proposition\!}

\def\proof{\noindent\emph{Proof}\ }
\def\endproof{\hfill$\square$}

\begin{document}

\title{\textbf{Gorenstein projective objects in homotopy categories}
\footnotetext{\noindent* Corresponding author.
\\\ \ \indent ~~ E-mail: lubo55@126.com (Lu B.), liuzk@nwnu.edu.cn (Liu Z.K.).\\\ \ \indent~~This research was
supported by National Natural Science Foundation of China
(No.11501451), Fundamental Research Funds for the Central Universities of China (No.31920150038) and XBMUYJRC (No.201406).}}
\author{\textbf{Lu Bo$^{\mathrm{a},*}$, Liu Zhongkui$^{\mathrm{b}}$} \\
\small $^a$College of Mathematics and Computer Science, Northwest University for Nationalities\\
\small Lanzhou 730124, Gansu, China\\
\small $^b$Department of Mathematics, Northwest Normal University\\
\small Lanzhou 730070, Gansu, China}
\date{}
\maketitle

\begin{abstract}

\vspace*{3mm}In this paper, we are concerned with Gorenstein projective objects in homotopy categories. Specifically, we present a characterization on Gorenstein projective objects in the category of complexes. Using this result, it is proved that the category of Gorenstein projective objects is a compactly generated triangulated category and an equivalence of triangulated categories is also given over some reasonably nice rings.

\vspace*{2mm} \noindent \emph{Key Words:} Gorenstein projective objects; Strongly Cartan-Eilenberg Gorenstein projective complex; Gorenstein projective module; Compactly generated triangulated category.

\vspace*{3mm} \noindent\emph{Mathematics Subject Classification:}
16E05; 18G05; 18E30; 18E10.
\end{abstract}

\vspace*{3mm} \setcounter{section}{1}
\subsection*{1. Introduction}

Triangulated categories were introduced by Grothendieck and Verdier in the early sixties
as the proper framework for doing homological algebra in an abelian category. Since then triangulated categories have found important applications in algebraic geometry, stable homotopy theory and homological
algebra.

Let $\mathcal{C}$ be a triangulated category with triangulation $\Delta$. Beligiannis [\cite{Be}] developed a
homological algebra in $\mathcal{C}$ which parallels the homological algebra in an exact category in
the sense of Quillen. He did this by specifying a class of triangles $\xi\subseteq\Delta$ which is closed
under translations and satisfies the analogous formal properties of a proper class of short
exact sequences. Such a class of triangles is called a proper class of triangles. By fixing
a proper class of triangles $\xi$, he defined projective and injective objects in triangulated categories.

In [\cite{AS}], Asadollahi and Salarian introduced and studied Gorenstein projective
and Gorenstein injective objects in triangulated categories. They are defined by modifying what Enochs and Jenda
have done in abelian categories [\cite{En95}]. These objects are called $\xi$-$\mathcal{G}$projective objects and $\xi$-$\mathcal{G}$injective objects, respectively.

A triangulated category $\mathcal{C}$ is said to have $enough~\xi$-$projectives$ if for any object $A \in \mathcal{C}$ there exists a triangle
$K \rightarrow P \rightarrow A\rightarrow\Sigma K$ in $\xi$ with $P$ $\xi$-projective. Dually one defines when $\mathcal{C}$ has $enough~\xi$-$injectives$.
In general it is not so easy to find a proper class $\xi$ of triangles in a triangulated category
having enough $\xi$-projectives or $\xi$-injectives.

Let $R$ be a ring, $\mathcal{K}(R)$ be the homotopy category of complexes of $R$-modules. Then $\mathcal{K}(R)$ is a triangulated category and so-called C-E projective complexes (or homotopically projective complexes) form the relative projective objects for a proper class of triangles in $\mathcal{K}(R)$, see [\cite{Be}, Section 12.4 and Section 12.5]. Motivated by the above, in the present article, we introduce and study the concept of Gorenstein projective objects using C-E projective complexes, and consider the relationship between Gorenstein projective objects and $\xi$-$\mathcal{G}$projective objects in [\cite{AS}]. Specifically, we get the following result as our main result in this note (cf. Theorem 3.3 and Theorem 3.12).

\begin{thm} An object $X$ in $\mathcal{K}(R)$ is~Gorenstein projective if and only if~$X$ is homotopy equivalent to~$G$, where $G_n$ is projective, $\mathrm{Z}_n(G)$, $\mathrm{B}_n(G)$ and~$\mathrm{H}_n(G)$ is~Gorenstein projective in~$R$-$\mathrm{Mod}$~for each~$n\in\mathbb{Z}$.
\end{thm}

\begin{thm} Let $G$ be an object in $\mathcal{K}(R)$. If $G$ is~Gorenstein projective object in Definition 3.2, then G is $\xi(\mathcal{P}_{\mathrm{CE}})$-$\mathcal{G}$projective in $\mathcal{K}(R)$.

\end{thm}

Compactly generatedness of triangulated categories is very important and useful. For instance, they allow the use of the Brown Representability Theorem and the Thomason Localization Theorem, both proved by Neeman in [\cite{N96}]. There are also many compactly generated triangulated categories. For example:

The stable module category $\underline{\mathrm{Mod}}(kG)$ of the group algebra of a finite
group over a field $k$; the compact objects are the objects induced by the finitely generated $kG$-modules, see [\cite{B07}].

The unbounded derived category of quasi-coherent sheaves over a quasi-compact separated scheme; the compact objects are the complexes of the
thick subcategory generated by the suspensions of powers of an ample line bundle, see [\cite{N96}].

The stable homotopy category of spectra; the compact objects are the finite CW-complexes, see [\cite{M83}].

The derived category $\mathcal{D}(R)$ of the abelian category $R$-Mod is
always compactly generated by the set $\mathcal{G} = \{S^n(R)| n\in\mathbb{Z}\}$, see [\cite{B07}].

However, surprisingly, the corresponding homotopy category $\mathcal{K}(R)$ is not even compactly generated when $R = \mathbb{Z}$, see [\cite{N96}, Lemma E.3.2]. J${\o}$rgensen proved in [\cite{J}] that the homotopy category of complexes of projective $R$-modules is compactly generated when the ring $R$ satisfies some
hypotheses. Neeman improved this result and showed that the homotopy category of complexes of projective $R$-modules is compactly generated when the ring $R$ is right coherent [\cite{N08}, Proposition 7.14]. In [\cite{K}] it is shown that the homotopy category of complexes of injective $R$-modules is compactly generated when $R$ is left noetherian, and this result is applied to give a new and interesting characterization of Gorenstein rings in terms of (totally) acyclic complexes of injective modules in [\cite{IK}, Corollary 5.5].

As an application of Theorem 1.1, we show that the category of all Gorenstein projective objects in $\mathcal{K}(R)$ is compactly generated over some reasonably nice rings, see Proposition 4.2. We also establish an equivalence of triangulated categories using Theorem 1.1, see Corollary 4.4.

To this end, we introduce and investigate the notion of strongly Cartan-Eilenberg Gorenstein projective complexes, see Definition 3.5.

It is an important question to establish relationships between a complex $X$ and the modules $X_n, \mathrm{Z}_n(X), n\in\mathbb{Z}$. It is well known that a complex $C$ is projective (respectively injective) if and only if $C$ is exact and $\mathrm{Z}_n(C)$ is projective (respectively injective) in $R$-Mod for each $n\in\mathbb{Z}$. In [\cite{En98}], Enochs and Garc\'{\i}a Rozas introduced and investigated the notion of Gorenstein projective and injective complexes. It was shown that a complex $C$ is Gorenstein projective (respectively Gorenstein injective) if and only if $C_m$ is Gorenstein projective (respectively Gorenstein injective) in $R$-$\mathrm{Mod}$ for $m\in\mathbb{Z}$ over $n$-Gorenstein rings. Liu and Zhang proved that Gorenstein injective version of this result holds over left Noetherian rings [\cite{Liu09}]. This has been further developed by Yang and Liu over any associated ring [\cite{Ya11}]. The following result is obtained in section 5 of this paper (cf. Proposition 3.8).

\begin{thm}
Let $G$ be a complex. Then the following statements are equivalent$:$

\noindent$(1)$ $G$ is strongly Cartan-Eilenberg Gorenstein projective$.$

\noindent$(2)$ $G_n$ is projective, $\mathrm{Z}_n(G)$, $\mathrm{B}_n(G)$ and $\mathrm{H}_n(G)$ are Gorenstein projective in R-$\mathrm{Mod}$ for each $n\in\mathbb{Z}.$

\end{thm}

We will use Foxby and Avramov's terminology [\cite{AV}] and a complex $C$ is $DG$-$projective$ (also see [\cite{Ga99}, p. 33]) if $C_n$ is projective and $\mathrm{Hom}_R(C,E)$ is exact for all exact complexes $E$. They are also termed semiprojective complexes. We get the following observation for strongly Cartan-Eilenberg Gorenstein projective complexes (cf. Proposition 5.5).

\begin{pro} Let G be an exact complex. Then the following statements are equivalent$:$

\noindent$(1)$ G is strongly Cartan-Eilenberg Gorenstein projective$.$

\noindent$(2)$ $G_n$ is projective in R-$\mathrm{Mod}$ for $n\in\mathbb{Z}$ and $\mathrm{Hom}_R(G,P)$ is exact for every Cartan-Eilenberg projective complex P$.$

\end{pro}

Note that projective complexes are strongly Cartan-Eilenberg Gorenstein projective complexes which are $\sharp$-projective by the definition of projective complexes and Theorem 1.3. However, the converse is not true in general, see Example 5.9 and 5.10. Recall that a complex $C$ is $\sharp$-$projective$ if $C_n$ is projective in $R$-Mod for each $n\in\mathbb{Z}$ [\cite{AV}].

%
%
%
%
%
%
%
%

As we know, an object $Q\in\mathcal{C}(R)$ is projective if and only if $Q$ is exact and $\mathrm{Z}_n(Q)$ is projective in $R$-$\mathrm{Mod}$ for each $n\in\mathbb{Z}$, where $\mathcal{C}(R)$ denotes the category of complexes of left $R$-modules. In [\cite{N08}], Neeman showed that the homotopy category of complexes of projective $R$-modules $\mathcal{K}(R$-Proj) is compactly generated when the ring $R$ is right coherent.
As an application of Theorem 1.4, we obtain a Gorenstein version of the above result in the corresponding category of complexes $\mathcal{C}(R$-Proj):

\begin{cor} Let $G$ be in $\mathcal{C}(R$-$\mathrm{Proj})$. Then $G$ is exact and $\mathrm{Z}_n(G)$ is Gorenstein projective in $R$-$\mathrm{Mod}$ for each $n\in\mathbb{Z}$ if and only if $\mathrm{Hom}_R(P,G)$ and $\mathrm{Hom}_R(G,P)$ are exact for every Cartan-Eilenberg projective complex $P.$

\end{cor}

Note that the above result can be considered as a new characterization of $\mathcal{A}$ complexes (see [\cite{Gi04}, Definition 3.3]) in $\mathcal{C}(R$-Proj) when $\mathcal{A}$ is the class of Gorenstein projective modules.

For the rest of the paper we will use the abbreviation C-E for Cartan-Eilenberg.

\stepcounter{section}
\subsection*{2. Preliminaries}

\vspace{3mm}Throughout this paper, $R$ denotes a ring with unity. A complex
\begin{center}$\cdots\stackrel{\delta_{2}}\longrightarrow
C_{1}\stackrel{\delta_{1}}\longrightarrow
C_{0}\stackrel{\delta_{0}}\longrightarrow
C_{-1}\stackrel{\delta_{-1}}\longrightarrow\cdots$
\end{center} of left $R$-modules will be denoted by $(C, \delta)$ or $C$. For such a complex $C$, we write
$\Sigma C$ for its suspension: $(\Sigma C)_n = C_{n-1}$ and $\delta^{\Sigma C} = -\delta^ C$. For a ring $R$, $R$-Mod denotes the category of
left $R$-modules.

We will use superscripts to distinguish complexes. So if
$\{C^i\}_{i\in I}$ is a family of complexes, $C^i$ will be
\begin{center}$\cdots\stackrel{\delta_{2}}\longrightarrow
C^i_{1}\stackrel{\delta_{1}}\longrightarrow
C^i_{0}\stackrel{\delta_{0}}\longrightarrow
C^i_{-1}\stackrel{\delta_{-1}}\longrightarrow\cdots.$
\end{center}

Given a left $R$-module $M$, we use the notation $D^m(M)$ to
denote the complex
\begin{center}$\cdots\stackrel{}\longrightarrow
0\stackrel{}\longrightarrow M\stackrel{id}\longrightarrow
M\stackrel{}\longrightarrow0\longrightarrow\cdots$
\end{center} with $M$ in the $m$th and $(m-1)$th positions. We also use the notation $S^m(M)$ to denote the complex with $M$ in
the $m$th place and 0 in the other places.

Given a complex $C$ and an integer $l$, the $l$th homology module of $C$
is the module $\mathrm{H}_l(C) = \mathrm{Z}_l(C)/\mathrm{B}_l(C)$ where $\mathrm{Z}_l(C)$ =
Ker$(\delta^C_l )$ and $\mathrm{B}_l(C)$ = Im$(\delta^C_{l+1})$.

Let $C$ and $D$ be complexes of left $R$-modules. We will denote by
Hom$_R(C,D)$ the complex of abelian groups with Hom$_R(C,D)_n=\prod\limits_{t\in\mathbb{Z}}\mathrm{Hom}_R(C_t,D_{n+t})$ and such that if $f\in \mathrm{Hom}_R(C,D)_n$ then $(\delta_n(f))_m=\delta^D_{m+n}f_m-(-1)^nf_{m+1}\delta^C_m$. $f$ is called a $chain$ $map$ of degree $n$ if $\delta_n(f)=0$. A chain map of degree 0 is called a $morphism$. We will use $\mathrm{Hom}(C,D)$ to denote the abelian
group of morphisms from $C$ to $D$ and Ext$^i$ for $i\geq0$ will
denote the groups we get from the right derived functor of Hom.

General background about complexes of $R$-modules can be found in [\cite{Ga99}].

Let $X$ and $Y$ be in $\mathcal{K}(R)$. We use $\mathrm{Hom}_{\mathcal{K}(R)}(X,Y)$ to denote the abelian group of morphisms from $X$ to $Y$.

We recall some notions and results needed in the paper.

A complex $P$ is $projective$ if the functor Hom$(P,-)$ is exact.

\begin{defn}[[\mbox{\cite{En95}, Definition 3.1]}] A left R-module M is called Gorenstein projective if there exists an exact sequence of projective left R-modules
$$\cdots\rightarrow P_{-1}\rightarrow
P_0\rightarrow P_1\rightarrow \cdots$$ with
$M\cong \mathrm{Ker}(P_0\rightarrow P_1)$ and which remains exact after applying $\mathrm{Hom}_R(-,P)$ for any projective left R-module P.
\end{defn}

\begin{defn}[\mbox{[\cite{Ve77}, p. 227]}] A complex $P$ is said to be $\mathrm{C}$-$\mathrm{E}$ projective if $P,
\mathrm{Z}(P), \mathrm{B}(P)$ and $\mathrm{H}(P)$ are complexes consisting of projective modules. Dually, ones defines the concept of $\mathrm{C}$-$\mathrm{E}$ injective complexes.

\end{defn}

Note that these complexes have origin in [\cite{CE56}] to give
the definitions of projective and injective resolutions of a complex of modules.

\begin{defn}[[\mbox{\cite{En11}, Definition 5.3]}] A sequence of complexes
$$\cdots\rightarrow C^{-1}\rightarrow C^{0}\rightarrow C^{1}\rightarrow\cdots$$
is said to be $\mathrm{C}$-$\mathrm{E}$ exact if

\noindent$(1)$ $\cdots\rightarrow
C^{-1}\rightarrow C^{0}\rightarrow C^{1}\rightarrow
\cdots$,

\noindent$(2)$ $\cdots\rightarrow
\mathrm{Z}(C^{-1})\rightarrow \mathrm{Z}(C^{0})\rightarrow \mathrm{Z}(C^{1})\rightarrow
\cdots$,

\noindent$(3)$ $\cdots\rightarrow
\mathrm{B}(C^{-1})\rightarrow \mathrm{B}(C^{0})\rightarrow \mathrm{B}(C^{1})\rightarrow
\cdots$,

\noindent$(4)$ $\cdots\rightarrow
C^{-1}/\mathrm{Z}(C^{-1})\rightarrow C^{0}/\mathrm{Z}(C^{0})\rightarrow
C^{1}/\mathrm{Z}(C^{1})\rightarrow \cdots$,

\noindent$(5)$ $\cdots\rightarrow
C^{-1}/\mathrm{B}(C^{-1})\rightarrow C^{0}/\mathrm{B}(C^{0})\rightarrow
C^{1}/\mathrm{B}(C^{1})\rightarrow \cdots$,

\noindent$(6)$ $\cdots\rightarrow
\mathrm{H}(C^{-1})\rightarrow \mathrm{H}(C^{0})\rightarrow \mathrm{H}(C^{1})\rightarrow
\cdots$

\noindent are all exact.
\end{defn}

%
%
%
%
%

In the following, we recall some definitions and results of triangulated categories used in this paper. The basic references for triangulated and derived categories are the original article of Verdier [\cite{Ve77}] and Hartshorne's notes [\cite{Ha66}]. Also Bernstein, Beilinson and Deligne [\cite{BE}] and Iversen [\cite{IV}] give introductions to these concepts. For terminology we shall follow [\cite{Be}].

Throughout this section we fix a triangulated category $\mathcal{C}=(\mathcal{C},T,\Delta)$, where $\mathcal{C}$ is an additive category, $T$ is the suspension functor, i.e., an autoequivalence of $\mathcal{C}$, and $\Delta$ is the triangulation.

A triangle (T) :$A\stackrel{f} \rightarrow B\stackrel{g} \rightarrow C\stackrel{ h} \rightarrow T A$ is called split if it is isomorphic to the triangle
$$A\stackrel{ \varepsilon} \longrightarrow A\oplus C\stackrel{\pi} \longrightarrow C\stackrel{ 0} \longrightarrow T A,$$
where $\varepsilon=\left(
                               \begin{array}{cc}
                                 1 & 0\\
                               \end{array}
                             \right),
\pi=\left(
      \begin{array}{c}
        0 \\
        1 \\
      \end{array}
    \right)
$.
It is easy to see that the above triangle is split if and only if $f$ has a retraction or $g$ has a section or
$h = 0$. The full subcategory of $\Delta$ consisting of the split triangles will be denoted by $\Delta_0$.

\begin{defn} [[\mbox{\cite{Be}, Definition 4.1]}] A full subcategory $\xi \subseteq Diag(\mathcal{C},T)$ is called a proper class of triangles if
the following conditions hold$:$

\noindent$(1)$ $\xi$ is closed under isomorphisms, finite coproducts and $\Delta_0 \subseteq \xi \subseteq \Delta$.

\noindent$(2)$ $\xi$ is closed under suspensions and is saturated.

\noindent$(3)$ $\xi$ is closed under base and cobase change.

\end{defn}

The definitions of base change, cobase change and saturated can be seen in [\cite{Be}, 2.1 and 2.2].

\begin{defn} [[\mbox{\cite{Be},  Definition 4.1]}]An object $P \in \mathcal{C}$ $($respectively $I \in \mathcal{C})$ is called $\xi$-projective $($respectively $\xi$-injective$)$ if for any triangle $A \rightarrow B \rightarrow C \rightarrow T A$ in $\xi$, the induced sequence
$$0\rightarrow \mathcal{C}(P,A)\rightarrow \mathcal{C}(P,B)\rightarrow \mathcal{C}(P,C)\rightarrow 0$$
$$(respectively~0\rightarrow \mathcal{C}(C, I )\rightarrow \mathcal{C}(B, I )\rightarrow \mathcal{C}(A, I )\rightarrow 0)$$
is exact in the category $\mathcal{A}b$ of abelian groups.
\end{defn}

In what follows, $\mathcal{P}(\xi$) (respectively $\mathcal{I}(\xi$)) will denote the full subcategory of $\xi$-projective
(respectively $\xi$-injective) objects of $\mathcal{C}$. It follows easily from the definition that the categories $\mathcal{P}(\xi$) and $\mathcal{I}(\xi$) are full, additive, closed under isomorphisms, direct summands and
$T$-stable.

\begin{defn} [[\mbox{\cite{AS},  Definition 3.1]}]An $\xi$-exact complex X is a diagram
$$\cdots\rightarrow X_1
\stackrel{d_1}\rightarrow X_0\stackrel{d_0}\rightarrow X_{-1}
\rightarrow\cdots$$
in $\mathcal{C}$, such that for integer n, there exist triangles
$$K_{n+1}
\stackrel{g_n}\rightarrow X_n
\stackrel{f_n}\rightarrow K_n
\stackrel{h_n}\rightarrow T K_{n+1}$$
in $\xi$ and a differential is defined as $d_n=g_{n-1}f_n$ for any n.
\end{defn}

\begin{defn} [[\mbox{\cite{AS},  Definition 3.2]}]A triangle $A\rightarrow B \rightarrow C\rightarrow T A$ in $\xi$, is called $\mathcal{C}(-,\mathcal{P}(\xi))$ exact, if for any
$Q\in \mathcal{P}(\xi)$, the induced complex
$$0\rightarrow \mathcal{C}(C,Q)\rightarrow \mathcal{C}(B,Q)\rightarrow \mathcal{C}(A,Q)\rightarrow0$$
is exact in $\mathcal{A}b$.
\end{defn}

\begin{defn} [[\mbox{\cite{AS},  Definition 3.3]}]A complete $\xi$-exact complex X is a diagram
$$X:\cdots\rightarrow X_1
\stackrel{d_1}\rightarrow X_0
\stackrel{d_0}\rightarrow X_{-1}\rightarrow\cdots$$
in $\mathcal{C}$ such that for all integers n, there exist $\mathcal{C}(-,P(\xi))$ exact triangles
$$K_{n+1}
\stackrel{g_n}\rightarrow X_n
\stackrel{f_n}\rightarrow K_n
\stackrel{h_n}\rightarrow T K_{n+1}$$
in $\xi$ where a differential $d_n$, for any integer n, is defined as $d_n=g_{n-1}f_n$.
\end{defn}

\begin{defn} [[\mbox{\cite{AS},  Definition 3.4]}] A complete $\xi$-projective resolution is a complete $\xi$-exact complex
$$P: \cdots\rightarrow P_1
\stackrel{d_1}\rightarrow P_0
\stackrel{d_0}\rightarrow P_{-1}\rightarrow\cdots$$
in $\mathcal{C}$ such that $P_n$, for any integer $n$, is $\xi$-projective.
\end{defn}

\begin{defn} [[\mbox{\cite{AS},  Definition 3.6]}]Let P be a complete $\xi$-projective resolution in $\mathcal{C}$. So for any integer n,
there exists a triangle
$$K_{n+1}
\stackrel{g_n}\rightarrow P_n
\stackrel{f_n}\rightarrow K_n
\stackrel{h_n}\rightarrow T K_{n+1}$$
in $\xi$. The objects $K_n$ for any integer n, are called $\xi$-$\mathcal{G}projective$.
\end{defn}

We denote by $\mathcal{GP}(\xi)$ the full subcategory of $\xi$-$\mathcal{G}$projective objects of $\mathcal{C}$. It follows
directly from the definition that the category $\mathcal{GP}(\xi)$ is full, additive and closed under isomorphisms. Every $\xi$-projective object is $\xi$-$\mathcal{G}$projective. In particular, there is an inclusion
of categories $\mathcal{P}(\xi)\subseteq \mathcal{GP}(\xi)$.

\stepcounter{section}

\subsection*{3. Gorenstein projective objects in homotopy categories}

In this section, we will introduce and study Gorenstein projective objects in homotopy categories. Specifically, we will consider the relationships between Gorenstein projective objects in $\mathcal{K}(R)$ and the corresponding objects in $\mathcal{C}(R)$.

Recall that an exact sequence of complexes $0\rightarrow A\rightarrow B\rightarrow C\rightarrow0$ is said to be $degreewise~split$ if $0\rightarrow A_n\rightarrow B_n\rightarrow C_n\rightarrow0$ is split in $R$-Mod for each $n\in\mathbb{Z}$.
\begin{defn}A sequence of complexes
$$\cdots\rightarrow C^{-1}\rightarrow C^{0}\rightarrow C^{1}\rightarrow\cdots$$
is said to be strongly $\mathrm{C}$-$\mathrm{E}$ exact if

\noindent$(1)$ $\cdots\rightarrow C^{-1}\rightarrow C^{0}\rightarrow C^{1}\rightarrow\cdots$ is $\mathrm{C}$-$\mathrm{E}$ exact$;$ and

\noindent$(2)$ $\cdots\rightarrow C^{-1}\rightarrow C^{0}\rightarrow C^{1}\rightarrow\cdots$ is degreewise split.
\end{defn}

\begin{defn}  An object $X$ in $\mathcal{K}(R)$ is said to be Gorenstein projective if there exists a strongly $\mathrm{C}$-$\mathrm{E}$ exact sequence of $\mathrm{C}$-$\mathrm{E}$ projective complexes in $\mathcal{C}(R)$
$$\mathbb{P}:\cdots\rightarrow P^{-2}\rightarrow P^{-1}\rightarrow
P^0\rightarrow P^1\rightarrow P^2\rightarrow\cdots$$ such that~$X\cong \mathrm{Ker}(P^0\rightarrow P^1)$~in~$\mathcal{K}(R)$ and which remains exact in $\mathcal{C}(R)$ after applying~$\mathrm{Hom}_{\mathcal{K}(R)}(-,Q)$ for any~$\mathrm{C}$-$\mathrm{E}$ projective complex~$Q$.

\end{defn}

In the following, we will use $\mathcal{K}(\mathcal{GP})$ to denote the class of Gorenstein projective objects in $\mathcal{K}(R)$.

The next theorem will give a characterization of Gorenstein projective objects.

\begin{thm} An object $X$ in $\mathcal{K}(R)$ is~Gorenstein projective if and only if~$X$ is homotopy equivalent to~$G$, where $G_n$ is projective, $\mathrm{Z}_n(G)$, $\mathrm{B}_n(G)$ and~$\mathrm{H}_n(G)$ is~Gorenstein projective in~$R$-$\mathrm{Mod}$~for each~$n\in\mathbb{Z}$.
\end{thm}
\vspace*{3mm}
$\mathbf{Proof~of~Theorem~3.3}$

\vspace*{3mm}
We will divide the proof of Theorem 3.3 into three steps.

$\mathbf{Step~1}$. We give the following lemma.

\begin{lem} Let~$0\rightarrow K\stackrel{g}\rightarrow M\stackrel{f}\rightarrow L\rightarrow0$~be a sequence of complexes,~$n\in\mathbb{Z}$. Then~$0\rightarrow K\rightarrow M\rightarrow L\rightarrow0$~is~$\mathrm{C}$-$\mathrm{E}$ exact if any two sequences of the following sequences are exact$:$

\noindent$(1)$~$0\rightarrow K_n\rightarrow M_n\rightarrow L_n\rightarrow0$.

\noindent$(2)$~$0\rightarrow \mathrm{Z}_n(K)\rightarrow \mathrm{Z}_n(M)\rightarrow \mathrm{Z}_n(L)\rightarrow0$.

\noindent$(3)$~$0\rightarrow \mathrm{B}_n(K)\rightarrow \mathrm{B}_n(M)\rightarrow \mathrm{B}_n(L)\rightarrow0$ or~$0\rightarrow K_n/\mathrm{Z}_n(K)\rightarrow M_n/\mathrm{Z}_n(M)\rightarrow L_n/\mathrm{Z}_n(L)\rightarrow0$.

\noindent$(4)$~$0\rightarrow K_n/\mathrm{B}_n(K)\rightarrow M_n/\mathrm{B}_n(M)\rightarrow L_n/\mathrm{B}_n(L)\rightarrow0$.

\noindent$(5)$~$0\rightarrow \mathrm{H}_n(K)\rightarrow \mathrm{H}_n(M)\rightarrow \mathrm{H}_n(L)\rightarrow0$.

\end{lem}
\proof If~(1)(2) or~(1)(3) or~(1)(4) or~(2)(3) or~(2)(5) or~(3)(4) or~(3)(5) or~(4)(5) hold, the result follows easily.

Suppose that~$0\rightarrow K_n\rightarrow M_n\rightarrow L_n\rightarrow0$ and~$0\rightarrow \mathrm{H}_n(K)\rightarrow \mathrm{H}_n(M)\rightarrow \mathrm{H}_n(L)\rightarrow0$ are exact. Let~$x\in \mathrm{B}_n(L)$. Since we have an exact commutative diagram:$$\xymatrix{
  M_{n+1} \ar[d]_{\delta^M_{n+1}} \ar[r]^{f_{n+1}} & L_{n+1} \ar[d]_{\delta^L_{n+1}} \ar[r]^{} & 0  \\
  M_{n}  \ar[r]^{f_{n}} & L_{n}  \ar[r]^{} & 0   ,}$$
there exists~$z\in L_{n+1}$ such that~$\delta_{n+1}^L(z)=x$. From~$f_{n+1}$~is epimorphic, there exists~$m\in M_{n+1}$ such that~$f_{n+1}(m)=z$. Thus ~$\delta_{n+1}^M(m)\in\mathrm{B}_n(M)$ and~$f_n(\delta_{n+1}^M(m))=\delta_{n+1}^Lf_{n+1}(m)=\delta_{n+1}^L(z)=x$, i.e., $f_n\mid_{\mathrm{B}_n(M)}:\mathrm{B}_n(M)\rightarrow \mathrm{B}_n(L)$ is epimorphic. On the other hand, using the exact commutative diagram:
$$\xymatrix{
  0 \ar[r]^{} & \mathrm{B}_n(M)\ar[d]_{f_{n}} \ar[r]^{} & \mathrm{Z}_n(M) \ar[d]_{f_{n}} \ar[r]^{}   & \mathrm{H}_n(M) \ar[d]_{\overline{f_{n}}} \ar[r]^{} & 0  \\
  0 \ar[r]^{} & \mathrm{B}_n(L) \ar[r]^{f_{n+1}} & \mathrm{Z}_n(L) \ar[r]^{f_{n}} & \mathrm{H}_n(L)  \ar[r]^{} & 0   ,}$$we get~$\mathrm{Z}_n(M)\rightarrow \mathrm{Z}_n(L)$ is an epimorphism. Hence~$0\rightarrow K\rightarrow M\rightarrow L\rightarrow0$~is~C-E exact by [\cite{En11} Lemma 5.2].

Let~$0\rightarrow \mathrm{Z}_n(K)\rightarrow \mathrm{Z}_n(M)\rightarrow \mathrm{Z}_n(L)\rightarrow0$ and~$0\rightarrow K_n/\mathrm{B}_n(K)\rightarrow M_n/\mathrm{B}_n(M)\rightarrow L_n/\mathrm{B}_n(L)\rightarrow0$ are exact. Then we have an exact commutative diagram:
$$\xymatrix{
  0 \ar[r]^{} & K_n/\mathrm{B}_n(K)\ar[d]^{\delta_{n}^{K/\mathrm{B}(K)}} \ar[r]^{} & M_n/\mathrm{B}_n(M) \ar[d]^{\delta_{n}^{M/\mathrm{B}(M)}} \ar[r]^{}   & L_n/\mathrm{B}_n(L) \ar[d]^{\delta_{n}^{L/\mathrm{B}(L)}} \ar[r]^{} & 0  \\
  0 \ar[r]^{} & K_{n-1}/\mathrm{B}_{n-1}(K)\ar[d]^{} \ar[r]^{} & M_{n-1}/\mathrm{B}_{n-1}(M)\ar[d]^{} \ar[r]^{} & L_{n-1}/\mathrm{B}_{n-1}(L)\ar[d]^{}  \ar[r]^{} & 0 \\
     & 0  & 0  & 0   &   ,}$$
where~$\mathrm{Ker}(\delta_{n}^{K/\mathrm{B}(K)})=\mathrm{H}_n(K)$, ~$\mathrm{Ker}(\delta_{n}^{M/\mathrm{B}(M)})=\mathrm{H}_n(M)$,~$\mathrm{Ker}(\delta_{n}^{L/\mathrm{B}(L)})=\mathrm{H}_n(L)$. By Snake Lemma,~$0\rightarrow \mathrm{H}_n(K)\rightarrow \mathrm{H}_n(M)\rightarrow \mathrm{H}_n(L)\rightarrow0$ is exact. Moreover,~$0\rightarrow \mathrm{Z}_n(K)\rightarrow \mathrm{Z}_n(M)\rightarrow \mathrm{Z}_n(L)\rightarrow0$ is exact, which implies that $0\rightarrow K\rightarrow M\rightarrow L\rightarrow0$~is~C-E exact.
\endproof

$\mathbf{Step~2}$. For convenience, we intorduce the following definition.

\begin{defn} A complex G is said to be strongly $\mathrm{C}$-$\mathrm{E}$ Gorenstein projective
if there is a strongly $\mathrm{C}$-$\mathrm{E}$ exact sequence of $\mathrm{C}$-$\mathrm{E}$ projective complexes
$$\mathbb{P}:\cdots\rightarrow P^{-1}\rightarrow
P^0\rightarrow P^1\rightarrow \cdots$$ such that
$G\cong \mathrm{Ker}(P^0\rightarrow P^1)$ and the functor
$\mathrm{Hom}(-,Q)$ leaves $\mathbb{P}$ exact whenever $Q$ is $\mathrm{C}$-$\mathrm{E}$
projective. In this case, $\mathbb{P}$ is called a strongly complete $\mathrm{C}$-$\mathrm{E}$ projective resolution of G.
\end{defn}

$\mathbf{Step~3}$. The following result, Proposition 3.8, will be established, which is a key observation for obtaining our main result. To this end, we first present the following lemmas.

\begin{lem} [[\mbox{\cite{En11}, Proposition 3.4]}]
A complex P is $\mathrm{C}$-$\mathrm{E}$ projective if and only if $P=P'\bigoplus P''$ where $P'$ is a projective complex and $P''$ is a graded module $($i.e. $\delta^{P''}=0)$ such that $P''\in\mathcal{C}(R$-$\mathrm{Proj})$, where $P''\in\mathcal{C}(R$-$\mathrm{Proj})$ denotes the category of those complexes each of whose terms is a projective module.
\end{lem}

\begin{lem} [[\mbox{\cite{Gi04}, Lemma 3.1]}]
For any left $R$-module $M$ and any complex of left R-modules $X$, we have the following natural isomorphisms$:$

\noindent$(1)$  $\mathrm{Hom}(D^{n}(M), X) \cong \mathrm{Hom}_{R}(M, X_{n}).$

\noindent$(2)$  $\mathrm{Hom}(S^{n}(M), X) \cong \mathrm{Hom}_{R}(M, \mathrm{Z}_{n}(X)).$

\noindent$(3)$   $\mathrm{Hom}(X, D^{n}(M)) \cong \mathrm{Hom}_{R}(X_{n-1}, M).$

\noindent$(4)$   $\mathrm{Hom}(X, S^{n}(M)) \cong \mathrm{Hom}_{R}(X_{n}/\mathrm{B}_{n}(X), M)$.

\end{lem}

\begin{pro}
Let $G$ be a complex. Then the following statements are equivalent$:$

\noindent$(1)$ $G$ is strongly $\mathrm{C}$-$\mathrm{E}$ Gorenstein projective$.$

\noindent$(2)$ $G_n$ is projective, $G_n/\mathrm{Z}_n(G)$ and $\mathrm{H}_n(G)$ are Gorenstein projective in R-$\mathrm{Mod}$ for each $n\in\mathbb{Z}.$

\noindent$(3)$ $G_n$ is projective, $\mathrm{Z}_n(G)$, $\mathrm{B}_n(G)$ and $\mathrm{H}_n(G)$ are Gorenstein projective in R-$\mathrm{Mod}$ for each $n\in\mathbb{Z}.$

\noindent$(4)$ $G_n$ is projective, $\mathrm{Z}_n(G)$, $\mathrm{B}_n(G)$, $\mathrm{H}_n(G)$, $G_n/\mathrm{B}_n(G)$ and $G_n/\mathrm{Z}_n(G)$ are Gorenstein projective in R-$\mathrm{Mod}$ for each $n\in\mathbb{Z}.$

\end{pro}

\proof $(1)\Rightarrow(2)$. Let $G$ has a strongly complete C-E projective resolution
$$\mathbb{P}=:\cdots\rightarrow P^{-1}\rightarrow
P^0\rightarrow P^1\rightarrow
\cdots.$$ Then the exact
sequence
$$\cdots\rightarrow P_n^{-1}\rightarrow
P_n^0\rightarrow P_n^1\rightarrow
\cdots$$
is split for each $n\in\mathbb{Z}$, which yields $G_n$ is projective in $R$-Mod for each $n\in\mathbb{Z}$.
We also get an exact sequence
$$\cdots\rightarrow P_n^{-1}/\mathrm{Z}_n(P^{-1})\rightarrow
P_n^0/\mathrm{Z}_n(P^{0})\rightarrow P_n^1/\mathrm{Z}_n(P^{1})\rightarrow\cdots\eqno(\dag_1)$$
of projective modules with $G_n/\mathrm{Z}_n(G)=\mathrm{Ker}(P_n^0/\mathrm{Z}_n(P^{0})\rightarrow
P_n^1/\mathrm{Z}_n(P^{1}))$ for all $n\in\mathbb{Z}$. So we only need to show that $\mathrm{Hom}_R(-,Q)$ leaves the sequence $(\dag_1)$ exact when $Q$ is a projective module.  Consider the exact sequence of complexes
 $$0\longrightarrow \textrm{B}_{n}(\mathbb{P})\longrightarrow \mathbb{P}_{n}\longrightarrow
 \mathbb{P}_{n}/\textrm{B}_{n}(\mathbb{P})\longrightarrow 0,\eqno(\dag_2)$$
 where
  $$\textrm{B}_{n}(\mathbb{P}) = \cdots \longrightarrow \textrm{B}_{n}(P^{1}) \longrightarrow
 \textrm{B}_{n}(P^{0}) \longrightarrow \textrm{B}_{n}(P^{-1}) \longrightarrow \cdots.$$ The sequence $(\dag_2)$ is degreewise split exact. For any projective module $Q$, we have the following exact sequence of complexes of $\mathbb{Z}$-modules
 $$0\longrightarrow \textrm{Hom}_{R}(\mathbb{P}_{n}/\textrm{B}_{n}(\mathbb{P}), Q)\longrightarrow
 \textrm{Hom}_{R}(\mathbb{P}_{n}, Q)\longrightarrow  \textrm{Hom}_{R}(\textrm{B}_{n}(\mathbb{P}), Q)\longrightarrow 0.$$
Since $\mathrm{Hom}(X, D^{n+1}(Q))\cong \mathrm{Hom}_{R}(X_{n}, Q)$ by Lemma 3.7, applying Hom$(-, D^{n+1}(Q))$ to the sequence $\mathbb{P}$ yields $\textrm{Hom}_{R}(\mathbb{P}_{n}, Q)$ is exact. Also, by applying Hom$(-, S^{n}(Q))$ to the sequence $\mathbb{P}$,
we get that
Hom$_{R}(\mathbb{P}_{n}/\textrm{B}_{n}(\mathbb{P}), Q)$ is exact by an argument analogous to the above.
Thus $\textrm{Hom}_{R}(\textrm{B}_{n}(\mathbb{P}), Q)$ is exact.  Therefore, $G_n/\mathrm{Z}_n(G)$ is Gorenstein projective in $R$-Mod since
$\mathrm{B}_{n-1}(P^i)\cong P_n^i/\mathrm{Z}_n(P^i)$. Similarly, we can also show that
$\mathrm{H}_n(G)$ is Gorenstein projective in $R$-Mod for each $n\in\mathbb{Z}$.

$(2)\Rightarrow(1)$. Since $G_n/\mathrm{Z}_n(G)$ and $\mathrm{H}_n(G)$ are
Gorenstein projective in $R$-Mod for all $n\in\mathbb{Z}$, $G_n/\mathrm{Z}_n(G)$ and $\mathrm{H}_n(G)$ have complete
projective resolutions. Then $G_n$ has a complete projective
resolution since $0\rightarrow \mathrm{H}_n(G)\rightarrow
G_n/\mathrm{B}_n(G)\rightarrow G_n/\mathrm{Z}_n(G)\rightarrow 0$ and $0\rightarrow
\mathrm{B}_n(G)\rightarrow G_n\rightarrow G_n/\mathrm{B}_n(G)\rightarrow 0$ are exact.
It follows that we can construct a strongly C-E exact sequence of
complexes$$\cdots\rightarrow P^{-1}\rightarrow
P^0\rightarrow P^1\rightarrow \cdots.\eqno(\dag_3)$$ By
the construction, we have that for each $n\in\mathbb{Z}$, the sequence
$$\cdots\rightarrow
P_n^{-1}\rightarrow P_n^0\rightarrow P_n^1\rightarrow
\cdots \eqno(\dag_4)$$
is exact with each $P_n^i$ projective in $R$-Mod, $G_n=\mathrm{Ker}(P_n^0\rightarrow P_n^1)$
and the sequence
$$\cdots\rightarrow
P_n^{-1}/\mathrm{B}_n(P^{-1})\rightarrow P_n^0/\mathrm{B}_n(P^{0})\rightarrow
P_n^1/\mathrm{B}_n(P^{1})\rightarrow
\cdots$$
is also exact with $P_n^i/\mathrm{B}_n(P^{i})$ projective, $G_n/\mathrm{B}_n(G)=\mathrm{Ker}(P_n^0/\mathrm{B}_n(P^0)\rightarrow P_n^1/\mathrm{B}_n(P^1))$ and $G_n/\mathrm{B}_n(G)$ is Gorenstein projective.
Since $G_n$ is projective for each $n\in\mathbb{Z}$, the exact sequence $(\dag_4)$ is split. Now we only show that the functor Hom$(-,P)$ leaves the strongly C-E exact
sequence $(\dag_3)$ exact when $P$ is C-E projective. Let $Q$ be a projective
module. Then \begin{center}$\cdots\rightarrow
\mathrm{Hom}(P_n^{1}/\mathrm{B}_n(P^{1}),Q)\rightarrow
\mathrm{Hom}(P_n^0/\mathrm{B}_n(P^{0}),Q)$

$\rightarrow\mathrm{Hom}(P_n^{-1}/\mathrm{B}_n(P^{-1}),Q)\rightarrow
\cdots$
\end{center}
is exact. Applying the functor Hom$(-,S^n(Q))$ to
the strongly C-E exact sequence $(\dag_3)$ and using Lemma 3.7 yield the exact sequence
$$\cdots\rightarrow
\mathrm{Hom}(P^1,S^n(Q))\rightarrow
\mathrm{Hom}(P^{0},S^n(Q))\rightarrow
\mathrm{Hom}(P^{-1},S^n(Q))\rightarrow\cdots.$$

We also apply the functor Hom$(-,D^{n+1}(Q))$ to the strongly C-E exact
sequence $(\dag_3)$ yields the sequence
$$\cdots\rightarrow
\mathrm{Hom}(P^1,D^n(Q))\rightarrow
\mathrm{Hom}(P^{0},D^n(Q))\rightarrow
\mathrm{Hom}(P^{-1},D^n(Q))\rightarrow\cdots$$
is exact similar to the above argument.
Note that $P=\bigoplus\limits_{i\in\mathbb{Z}}(D^i(Q_i)\oplus S^i(Q'_i))$ for any C-E projective complex $P$ by Lemma 3.6, where $Q_i$ and $Q'_i$ are projective in $R$-Mod. Hence (2) follows.

$(2)\Leftrightarrow(3)\Leftrightarrow(4)$ follows by [\cite{Ho04},Theorem 2.5].
\endproof

\vspace{2mm}

Note that $H_n(\mathrm{Hom}_R(X,Y))\cong \mathrm{Hom}_{\mathcal{K}(R)}(X,\Sigma^{-n}Y)$. Thus Theorem 3.3 immediately from Lemma 3.4, Definition 3.5 and Proposition 3.8.\endproof

\vspace{2mm}

From [\cite{Be}, Section 12.4 and Section 12.5], C-E projective complexes (or homotopically projective complexes) form the relative projective objects for a proper class of triangles in $\mathcal{K}(R)$. It is easy to see that $C$ is C-E projective (or homotopically projective) in $\mathcal{K}(R)$ if and only if $C$ is C-E projective in $\mathcal{C}(R)$ (see, Definition 2.2). In the following, we will use $\mathcal{P}_{\mathrm{CE}}$ to denote the class of C-E projective complexes (or homotopically projective complexes) in $\mathcal{K}(R)$. We will use $\xi(\mathcal{P}_{\mathrm{CE}})$ to denote the proper class which corresponds to C-E projective complexes. Take $\xi=\xi(\mathcal{P}_{\mathrm{CE}})$. Then $\xi$-$\mathcal{G}$projective objects in [\cite{AS}]
will be said by $\xi(\mathcal{P}_{\mathrm{CE}})$-$\mathcal{G}$projective.

In order to establish a relationship between Gorenstein projective objects and $\xi$-$\mathcal{G}$projective objects in [\cite{AS}], we first give the following results.

\begin{lem} Let $\mathcal{K}(R)$ be a homotopy category of complexes. Then a triangle $D\rightarrow F\rightarrow C\rightarrow T D$ is in $\xi(\mathcal{P}_{\mathrm{CE}})$ if and only if $0\rightarrow \mathrm{Hom}_{\mathcal{K}(R)}(P,D)\rightarrow \mathrm{Hom}_{\mathcal{K}(R)}(P,F)\rightarrow \mathrm{Hom}_{\mathcal{K}(R)}(P,C)\rightarrow 0$ is exact for any $P\in\mathcal{P}_{\mathrm{CE}}$.

\end{lem}
\proof It follows from [\cite{Be}, Lemma 4.2].\endproof

\begin{lem} Let $\mathcal{K}(R)$ be a homotopy category of complexes. Then $$\xi(\mathcal{P}_{\mathrm{CE}})=\{X\rightarrow Y\rightarrow Z\rightarrow T X~|~X\rightarrow Y\rightarrow Z\rightarrow T X\cong A\rightarrow B\rightarrow C\rightarrow T A~in\ \mathcal{K}(R)~\},$$
where~$~0\rightarrow A\rightarrow B\rightarrow C\rightarrow0$~is strongly ~$\mathrm{C}\textrm{-}\mathrm{E}$~exact~in~$\mathcal{C}(R)$.

\end{lem}
\proof Take $\eta=\{X\rightarrow Y\rightarrow Z\rightarrow T X~|~X\rightarrow Y\rightarrow Z\rightarrow T X\cong A\rightarrow B\rightarrow C\rightarrow T A~in\ \mathcal{K}(R)~\},$
where~$~0\rightarrow A\rightarrow B\rightarrow C\rightarrow0$~is strongly $\mathrm{C}\textrm{-}\mathrm{E}$~exact~in~$\mathcal{C}(R)$.

Let $M\rightarrow N\rightarrow L\rightarrow T M$ be in $\xi(\mathcal{P}_{\mathrm{CE}})$. Then there exists a degreewise split exact sequence $0\rightarrow R\rightarrow S\rightarrow P\rightarrow0$ such that $M\rightarrow N\rightarrow L\rightarrow T M$ is isomorphic to
$ R\rightarrow S\rightarrow P\rightarrow TR$ in $\mathcal{K}(R)$.

On the other hand,$$0\rightarrow \mathrm{Hom}_{\mathcal{K}(R)}(Q,\Sigma^{-n}R)\rightarrow \mathrm{Hom}_{\mathcal{K}(R)}(Q,\Sigma^{-n}S)\rightarrow \mathrm{Hom}_{\mathcal{K}(R)}(Q,\Sigma^{-n}P)\rightarrow0$$ is exact for any C-E projective complex $Q$, and so$$0\rightarrow \mathrm{H}_n(\mathrm{Hom}_R(Q,R))\rightarrow \mathrm{H}_n(\mathrm{Hom}_R(Q,S))\rightarrow \mathrm{H}_n(\mathrm{Hom}_R(Q,P))\rightarrow0$$ is exact.
 Moreover, $\mathrm{S}^0(R)$ is C-E projective, which yields that $0\rightarrow \mathrm{H}_n(R)\rightarrow \mathrm{H}_n(S)\rightarrow \mathrm{H}_n(P)\rightarrow0$ is exact. Therefore, ~$\xi(\mathcal{P}_{\mathrm{CE}}) \subseteq~\eta$.

Let $0\rightarrow A\rightarrow B\rightarrow C\rightarrow0$~be a~strongly~$\mathrm{C}$-$\mathrm{E}$~exact~sequence of complexes. Then $0\rightarrow A_i\rightarrow B_i\rightarrow C_i\rightarrow0$ is split for $i\in\mathrm{Z}$. This means $A\rightarrow B\rightarrow C\rightarrow \Sigma A$ is a triangle in $\mathcal{K}(R)$.

Now suppose that for any $P\in\mathcal{P}_{\mathrm{CE}}$. There is a C-E projective complex $Q$ such that $P\cong Q$ in $\mathcal{K}(R)$. Then we get the following exact sequence of Abelian groups
$$0\rightarrow \mathrm{Hom}(Q,A)\rightarrow \mathrm{Hom}(Q,B)\rightarrow \mathrm{Hom}(Q,C)\rightarrow0.$$
Thus
$$0\rightarrow \mathrm{Z}_n(\mathrm{Hom}_R(Q,A))\rightarrow \mathrm{Z}_n(\mathrm{Hom}_R(Q,B))\rightarrow \mathrm{Z}_n(\mathrm{Hom}_R(Q,C))\rightarrow0$$
is exact.

Note that the sequence $0\rightarrow \mathrm{Hom}_R(Q,A)\rightarrow \mathrm{Hom}_R(Q,B)\rightarrow \mathrm{Hom}_R(Q,C)\rightarrow0$
is exact since $0\rightarrow A\rightarrow B\rightarrow C\rightarrow0$ is degreewise split. Then $0\rightarrow \mathrm{Hom}_R(Q,A)\rightarrow \mathrm{Hom}_R(Q,B)\rightarrow \mathrm{Hom}_R(Q,C)\rightarrow0$ is C-E exact. We get that
$$0\rightarrow \mathrm{H}_0(\mathrm{Hom}_R(Q,A))\rightarrow \mathrm{H}_0(\mathrm{Hom}_R(Q,B))\rightarrow \mathrm{H}_0(\mathrm{Hom}_R(Q,C))\rightarrow0$$ is exact, i.e., $0\rightarrow \mathrm{Hom}_{\mathcal{K}(R)}(Q,A)\rightarrow \mathrm{Hom}_{\mathcal{K}(R)}(Q,B)\rightarrow \mathrm{Hom}_{\mathcal{K}(R)}(Q,C)\rightarrow0$ is exact. Therefore, $$0\rightarrow \mathrm{Hom}_{\mathcal{K}(R)}(P,A)\rightarrow \mathrm{Hom}_{\mathcal{K}(R)}(P,B)\rightarrow \mathrm{Hom}_{\mathcal{K}(R)}(P,C)\rightarrow0$$ is exact, which concludes that~$\eta \subseteq\xi(\mathcal{P}_{\mathrm{CE}})$.
\endproof

\begin{cor} Let P be a complex. Then the following conditions are equivalent:

\noindent$(1)$~$P\in\mathcal{P}_{\mathrm{CE}}$, i.e.,~$P$ is homotopy equivalent to a complex having projective components and zero differential.

\noindent$(2)$ $P$ is homotopy equivalent to~$\mathrm{C}$-$\mathrm{E}$ projective complex.

\noindent$(3)$ $0\rightarrow \mathrm{Hom}_{\mathcal{K}(R)}(P,A)\rightarrow \mathrm{Hom}_{\mathcal{K}(R)}(P,B)\rightarrow \mathrm{Hom}_{\mathcal{K}(R)}(P,C)\rightarrow0$ is exact for every strongly~$\mathrm{C}$-$\mathrm{E}$ exact sequence~$0\rightarrow A\rightarrow B\rightarrow C\rightarrow0$.

\noindent$(4)$ $0\rightarrow \mathrm{Hom}_{\mathcal{C}(R)}(P,A)\rightarrow \mathrm{Hom}_{\mathcal{C}(R)}(P,B)\rightarrow \mathrm{Hom}_{\mathcal{C}(R)}(P,C)\rightarrow0$ is exact for every strongly~$\mathrm{C}$-$\mathrm{E}$ exact sequence~$0\rightarrow A\rightarrow B\rightarrow C\rightarrow0$.
\end{cor}
\proof (1)$\Rightarrow$(2) is obviously.

(2)$\Rightarrow$(1) Let $Q$ be a C-E projective. Then $Q=Q^1\oplus Q^2$, where $Q^1$ is a projective complex, $Q^2_i$ is a projective module and $\delta_i^{Q^2}=0$ for $i\in\mathbb{Z}$. Thus the exact sequence~$0\rightarrow Q^1\rightarrow Q\stackrel{f}\rightarrow Q^2\rightarrow0$ is degreewise split. Since~$Q^1$ is contractible,
$f:Q\rightarrow Q^2$ is a homotopy equivalence by~[\cite[Lemma~3.4.8]{CH11}].

(1)$\Leftrightarrow$(3) follows from Lemma 3.10.

(3)$\Leftrightarrow$(4) Since~$0\rightarrow A\rightarrow B\rightarrow C\rightarrow0$ is a degreewise split exact sequence, $0\rightarrow \mathrm{Hom}_{R}(P,A)\rightarrow \mathrm{Hom}_{R}(P,B)\rightarrow \mathrm{Hom}_{R}(P,C)\rightarrow0$ is exact. Hence the result follows.\endproof

\begin{thm} Let $G$ be an object in $\mathcal{K}(R)$. If $G$ is~ a Gorenstein projective object in Definition 3.2, then G is $\xi(\mathcal{P}_{\mathrm{CE}})$-$\mathcal{G}$projective in $\mathcal{K}(R)$.

\end{thm}
\proof Since $G$ is a Gorenstein projective object in $\mathcal{K}(R)$, there exists a strongly C-E Gorenstein projective complex $X$ by Theorem 3.3. Then there is a strongly C-E exact sequence
$$0\rightarrow K^1\rightarrow P^0\rightarrow X\rightarrow0$$
with $P^0$ C-E projective, $K^1$ strongly C-E Gorenstein projective in $\mathcal{C}(R)$ and $K^1\rightarrow P^0\rightarrow X\rightarrow\Sigma K^1$ is a triangle in
$\mathcal{K}(R)$. By the definition of strongly C-E Gorenstein projective complexes, we get that
$$0\rightarrow \mathrm{Hom}(X,Q)\rightarrow \mathrm{Hom}(P^0,Q)\rightarrow \mathrm{Hom}(K^1,Q)\rightarrow0$$
is exact for any C-E projective complex $Q$, which yields the sequence
$$0\rightarrow \mathrm{Z}_0(\mathrm{Hom}_R(X,Q))\rightarrow \mathrm{Z}_0(\mathrm{Hom}_R(P^0,Q))\rightarrow \mathrm{Z}_0(\mathrm{Hom}_R(K^1,Q))\rightarrow0$$
is exact. On the other hand,
$0\rightarrow \mathrm{Hom}_R(X,Q)\rightarrow \mathrm{Hom}_R(P^0,Q)\rightarrow \mathrm{Hom}_R(K^1,Q)\rightarrow0$ is exact since $0\rightarrow K^1\rightarrow P^0\rightarrow X\rightarrow0$ is degreewise split. Thus $$0\rightarrow \mathrm{Hom}_R(X,Q)\rightarrow \mathrm{Hom}_R(P^0,Q)\rightarrow \mathrm{Hom}_R(K^1,Q)\rightarrow0$$ is C-E exact, and so $$0\rightarrow \mathrm{H}_0(\mathrm{Hom}_R(G,Q))\rightarrow \mathrm{H}_0(\mathrm{Hom}_R(P^0,Q))\rightarrow \mathrm{H}_0(\mathrm{Hom}_R(K^1,Q))\rightarrow0$$ is exact, which means $0\rightarrow \mathrm{Hom}_{\mathcal{K}(R)}(X,Q)\rightarrow \mathrm{Hom}_{\mathcal{K}(R)}(P^0,Q)\rightarrow \mathrm{Hom}_{\mathcal{K}(R)}(K^1,Q)\rightarrow0$ is exact. We can also show that there is a triangle $X\rightarrow P^{-1}\rightarrow L^{-1}\rightarrow\Sigma X$ in $\mathcal{K}(R)$ and $0\rightarrow \mathrm{Hom}_{\mathcal{K}(R)}(L^{-1},Q)\rightarrow \mathrm{Hom}_{\mathcal{K}(R)}(P^{-1},Q)\rightarrow \mathrm{Hom}_{\mathcal{K}(R)}(X,Q)\rightarrow0$ is exact for any C-E projective complex $Q$. Therefore, $X$ is $\xi(\mathcal{P}_{\mathrm{CE}})$-$\mathcal{G}$projective in $\mathcal{K}(R)$. Note that $G\cong X$ in $\mathcal{K}(R)$. Thus $G$ is $\xi(\mathcal{P}_{\mathrm{CE}})$-$\mathcal{G}$projective in $\mathcal{K}(R)$.\endproof


As a direct consequence of Theorem 3.3 and Theorem 3.12, we get the following result.

\begin{cor}  Let $G$ be an object in $\mathcal{K}(R)$.

\noindent$(1)$  G is $\xi(\mathcal{P}_{\mathrm{CE}})$-$\mathcal{G}$projective in $\mathcal{K}(R)$.

\noindent$(2)$ $G_n$ is projective, $G_n/\mathrm{Z}_n(G)$ and $\mathrm{H}_n(G)$ are Gorenstein projective in R-$\mathrm{Mod}$ for each $n\in\mathbb{Z}.$

\noindent$(3)$ $G_n$ is projective, $\mathrm{Z}_n(G)$, $\mathrm{B}_n(G)$ and $\mathrm{H}_n(G)$ are Gorenstein projective in R-$\mathrm{Mod}$ for each $n\in\mathbb{Z}.$

\noindent$(4)$ $G_n$ is projective, $\mathrm{Z}_n(G)$, $\mathrm{B}_n(G)$, $\mathrm{H}_n(G)$, $G_n/\mathrm{B}_n(G)$ and $G_n/\mathrm{Z}_n(G)$ are Gorenstein projective in R-$\mathrm{Mod}$ for each $n\in\mathbb{Z}.$

Then $(2)\Longleftrightarrow(3)\Longleftrightarrow(4)\Longrightarrow(1)$.

\end{cor}

\begin{rem}In [\cite{AS}], Asadollahi and Salarian introduce and investigate $\xi$-$\mathcal{G}$projective objects in triangulated categories. From Corollary 3.13, we obtain a specific example of $\xi$-$\mathcal{G}$projective objects in $\mathcal{K}(R)$.
\end{rem}

\stepcounter{section}

\subsection*{4. Applications}

In this section, we will give some applications to our main result, Theorem 3.3, Theorem 3.12 in section 3.

\begin{defn}[[\mbox{\cite{N01}, Definition 2.1]}] Let $\mathcal{C}$ be a triangulated category, closed under set-indexed coproducts. An object $C\in \mathcal{C}$ is compact if the natural map$$\coprod\limits_{i\in I}\mathrm{Hom}_{\mathcal{C}}(C,X_i)\longrightarrow \mathrm{Hom}_{\mathcal{C}}(C,\coprod\limits_{i\in I}X_i)$$
is an isomorphism for any family $\{X_i\}_{i\in I}$ of objects in $\mathcal{C}$. A set of objects $\mathcal{G}\subseteq \mathcal{C}$ is
called a generating set if the implication
$$\mathrm{Hom}_{\mathcal{C}}(G,X) = 0~for~all~G\in\mathcal{G}\Longrightarrow  X\cong0$$
holds for all $X\in \mathcal{C}$. If $\mathcal{C}$ has a generating set consisting of compact objects, then $\mathcal{C}$
is called compactly generated.
\end{defn}

Recall that a ring $R$ is $0$-Gorenstein if $R$ is left and right Noetherican and id$(_RR)=0$.

\begin{pro} Let $R$ be a $0$-Gorenstein ring. Then the category $\mathcal{K}(\mathcal{GP})$ is a compactly generated triangulated category.
\end{pro}
\proof Since $R$ is a $0$-Gorenstein ring, it is easy to check that $\mathcal{K}(\mathcal{GP})$ is a triangulated category. Note that the category $\mathcal{K}(\mathcal{GP})$  happens to be the homotopy category of complexes of projective modules $\mathcal{K}(R$-Proj) over $0$-Gorenstein rings by Theorem 3.3. Using [\cite{N08}, Proposition 7.14], the result follows.
\endproof

We collect in the following some $0$-Gorenstein rings, see [\cite{R92}].
\begin{rem}

\noindent$(1)$  If G is a finite group and k is any field, then the group rings kG is $0-$Gorenstein.

\noindent$(2)$ If R is a PID and I is a nonzero proper ideal, then R/I is $0-$Gorenstein.

\noindent$(3)$ The rings $\mathbb{I}_n$ $(integers~mod~n)$, where $n>1$, and the rings k[x]/I, where k is a field and I nonzero ideal, are $0-$Gorenstein rings.
\end{rem}

Dualizing complexes are popular gadgets in homological algebra. Recall that a dualizing complex for $R$ is a complex $D$ of $R$-modules with the following properties:

(a) The cohomology of $D$ is bounded and finitely generated over $R$.

(b) The injective dimension id$_RD$ is finite.

(c) The canonical morphism $R\longrightarrow \mathrm{RHom}_R(D,D)$ in the derived category D$(R)$ is an isomorphism.

We write $K^{+,b}(R$-proj) for the subcategory of $\mathcal{K}(R)$ consisting of complexes $X$ of finitely generated projective modules with H$(X)$ bounded and $X_n = 0$ for $n\ll0$, and D$^f (R)$ for its image in D$(R)$, the derived category of $R$-modules.

\begin{cor} Let R be a noetherian ring with a dualizing complex D. Then

\noindent$(1)$ The subcategory of $\mathcal{E\mathcal{K}(\mathcal{GP})}$ is a compactly generated triangulated category, where $\mathcal{E\mathcal{K}(\mathcal{GP})}$ denotes the subcategory of all exact objects in $\mathcal{K}(\mathcal{GP})$.

\noindent$(2)$ There is an equivalence $$\mathrm{D}^f(R)/\mathrm{Thick}(R,D)\stackrel{\sim}\rightarrow{\mathcal{EK}^c(\mathcal{GP})}^{op},$$ where $\mathrm{Thick}(R,D)$ is the thick subcategory of $\mathrm{D}^f(R)$ generated by $R$ and $D$, $\mathcal{EK}^c(\mathcal{GP})$  is the class of compact objects in $\mathcal{E\mathcal{K}(\mathcal{GP})}$.
\end{cor}
\proof It follows by [\cite{IK}, Theorem 5.3] and Theorem 3.3.
\endproof

Let $\mathcal{X}$ be a class of objects in a triangulated category $\mathcal{C}$. Recall that the full subcategories $$\mathcal{X}^\perp = \{Y\in\mathcal{C} | \mathrm{Hom}_\mathcal{C}(T^nX, Y ) = 0~\mathrm{for~all}~X\in \mathcal{X}~\mathrm{and}~n\in\mathbb{Z}\}$$
is called the right orthogonal class to $\mathcal{X}$.

Let $\mathrm{i}R$ be an injective resolution of $R$ and $D^\ast = \mathrm{S}(\mathrm{i}R)$, where S = q$\circ\mathrm{Hom}_R(D,-) : \mathcal{\mathcal{K}}(R$-Inj)$\longrightarrow \mathcal{K}(R$-Proj), q denotes the right adjoint of the inclusion $\mathcal{K}(R$-$\mathrm{Proj})\longrightarrow \mathcal{K}(R$-Flat).

The following result gives a characterization of a special subcategory of $\mathcal{EK}^c(\mathcal{GP})$.

\begin{cor} Let R be a noetherian ring with a dualizing complex D. Then $$\mathcal{EK}^c(\mathcal{GP})=\{R,D^\ast\}^\perp.$$
\end{cor}
\proof It follows by [\cite{IK}, Proposition 5.9] and Theorem 3.3.
\endproof

We would like already here to point out that many rings have dualizing complexes, see [\cite{PJ,IK}].
\begin{rem}

\noindent$(1)$ A noetherian local commutative ring has a dualizing complex if and only
if it is a quotient of a Gorenstein noetherian local commutative ring.

\noindent$(2)$ Any local ring of a scheme of locally finite type over a field has a dualizing complex.

\noindent$(3)$ Any complete noetherian local commutative ring has a dualizing complex.

\noindent$(4)$ Let R be a commutative noetherian local ring with maximal ideal $\mathfrak{m}$ and residue
field $k = A/\mathfrak{m}$. Assume that $\mathfrak{m}^2 = 0$, and that $\mathrm{rank}_k(\mathfrak{m})\geq2$. Observe that R is not
Gorenstein$;$ for instance, its socle is $\mathfrak{m}$, and hence of rank at least $2$. Let E denote the
injective hull of the R-module k$;$ this is a dualizing complex for R.

\end{rem}

\stepcounter{section}

\subsection*{5. Some notes on strongly Cartan-Eilenberg Gorenstein projective complexes}

In this section, we will investigate further the notion of strongly C-E Gorenstein projective complexes which plays an important role in the study of Gorenstein projective objects in $\mathcal{K}(R)$.

We first establish several basic facts.

By [\cite{En11}, Proposition 6.3], we can compute derived functors of Hom$(-,-)$
using C-E projective resolutions or C-E injective resolutions. For given $C$ and $D$ we will
denote these derived functors applied to $(C, D)$ as
$\overline{\mathrm{Ext}}^n(C, D)$. It is obvious that
$\overline{\mathrm{Ext}}^1(C, D)\subseteq \mathrm{Ext}^1(C, D)$.

%

The first statement can be obtained by [\cite{En11}, Proposition 6.3], and the second follows from the definition of strongly C-E Gorenstein projective complexes and [\cite{En11}, Proposition 10.1].

\begin{lem}
\noindent$(1)$ A complex P is $\mathrm{C}$-$\mathrm{E}$ projective if and only if $\overline{\mathrm{Ext}}^i(P,D)=0$ for any complex D
and all $i\geq1.$

\noindent$(2)$ If P is a $\mathrm{C}$-$\mathrm{E}$ projective complex and G is a strongly $\mathrm{C}$-$\mathrm{E}$ Gorenstein projective complex, then $\overline{\mathrm{Ext}}^n(G,P)=0$ for each $n\geq 1$. Moreover, $\mathrm{Ext}^1(G,P)=0$ whenever G is exact.

\end{lem}

\begin{lem} Let G be a complex with $G_i$ Gorenstein projective in R-$\mathrm{Mod}$ for each $i\in\mathbb{Z}$. Then $\mathrm{Hom}_R(G,P)$ is exact for any $\mathrm{C}$-$\mathrm{E}$ projective complex P if and only if $\mathrm{Ext}^1(G,P)=0$ for any $\mathrm{C}$-$\mathrm{E}$ projective complex P.

\end{lem}
\proof It follows from [\cite{Gi04}, Lemma 2.1]. \endproof

\begin{lem}[[\mbox{\cite{GL}, Lemma 2.4]}]Let $0\rightarrow M\stackrel{f} \rightarrow P \rightarrow N \rightarrow 0$ be a short exact sequence of modules.
If M and P are projective and N is Gorenstein projective, then $\mathrm{Coker}(\alpha)$ is Gorenstein projective for
any homomorphism $f':M \rightarrow P'$ with $P'$ projective, where $\alpha= (f, f'): M \rightarrow P \oplus P'$ is defined by
$\alpha(x) = (f(x), f'(x))$ for any $x\in M.$

\end{lem}

\begin{lem} [[\mbox{\cite{LU14}, Lemma 3.2]}] Let $0\rightarrow A\rightarrow B \rightarrow C \rightarrow 0$ be a short exact sequence of complexes with A exact. Then $0\rightarrow A\rightarrow B \rightarrow C \rightarrow 0$ is $\mathrm{C}$-$\mathrm{E}$ exact.

\end{lem}

\vspace{3mm}Inspired by the notion of DG-projective complexes, see the introduction, we obtain the following result.

\begin{pro} Let G be an exact complex. Then the following statements are equivalent$:$

\noindent$(1)$ G is strongly $\mathrm{C}$-$\mathrm{E}$ Gorenstein projective$.$

\noindent$(2)$ $G_n$ is projective in R-$\mathrm{Mod}$ for each $n\in\mathbb{Z}$ and $\mathrm{Hom}_R(G,P)$ is exact for every $\mathrm{C}$-$\mathrm{E}$ projective complex P$.$

\end{pro}
\proof $(1)\Rightarrow(2)$. It follows from Proposition 3.8, Lemma 5.1 and Lemma 5.2.

$(2)\Rightarrow(1)$.  Since $G_n$ is projective, there exists an exact sequence of modules $$0\rightarrow G_n\stackrel{f_n}\rightarrow X_n\rightarrow H_n\rightarrow0$$
with $X_n$ projective and $H_n$ Gorenstein projective for each $n\in\mathbb{Z}$. Put
$$P^0=:\cdots\rightarrow {P_{n+1}^0}\stackrel{\delta_{n+1}^{P^0}} \rightarrow {P_{n}^0} \stackrel{\delta_{n}^{P^0}}\rightarrow {P_{n-1}^0}\stackrel{\delta_{n-1}^{P^0}}\rightarrow {P_{n-2}^0}\rightarrow\cdots$$ with ${P_n^0}=X_n\oplus X_{n-1}$ and $\delta_{n}^{P^0}:{P_{n}^0} \rightarrow {P_{n-1}^0}$ defined via $\delta_{n}^{P^0}(x,y)=(y,0)$ for $n\in\mathbb{Z}$ and any $(x,y)\in X_n\oplus X_{n-1}$. Clearly, $P^0$ is exact and Z$_n(P^0)$ is projective in $R$-Mod for each $n\in\mathbb{Z}$. Then $P^0$ is C-E projective. Now we have a morphism $\alpha=(\alpha_n)_{n\in\mathbb{Z}}:G \rightarrow P^0$ of complexes as follows:
$$\xymatrix{
  \cdots \ar[r]^{}& G_{n+1} \ar[d]_{(f_{n+1},f_n\delta^G_{n+1})} \ar[r]^{\delta^G_{n+1}} & G_n \ar[d]_{(f_{n},f_{n-1}\delta^G_{n})} \ar[r]^{\delta^G_{n}} & G_{n-1} \ar[d]_{(f_{n-1},f_{n-2}\delta^G_{n-1})} \ar[r]^{} &  \cdots \\
 \cdots \ar[r]^{} &  X_{n+1}\oplus X_{n} \ar[r]^{\delta_{n+1}^{P^0}} &  X_n\oplus X_{n-1} \ar[r]^{\delta_{n}^{P^0}} &  X_{n-1}\oplus X_{n-2} \ar[r]^{} & \cdots .  }
$$
It is clear that $\alpha$ is injective and so we have a short exact sequence of complexes
$$0\rightarrow G\stackrel{\alpha}\rightarrow P^0\rightarrow L^1\rightarrow 0,\eqno(\ddag_1)$$
where $L^1=\mathrm{Coker}(\alpha)$ and the sequence $(\ddag_1)$ is C-E exact by Lemma 5.4. Using Lemma 5.3, $L^1_n$
is Gorenstein projective for each $n\in\mathbb{Z}$, which means the exact sequence $0\rightarrow G_n\stackrel{\alpha_n}\rightarrow P_n^0\rightarrow L_n^1\rightarrow 0$ is split. And so we get that the sequence of abelian groups
$$0 \rightarrow \mathrm{Hom}_R(L^1,P) \rightarrow \mathrm{Hom}_R(P^0,P) \rightarrow \mathrm{Hom}_R(G,P) \rightarrow 0$$
is exact for any C-E projective complex $P$. Since Ext$^1(P^0,P)=0$ by Lemma 5.1 and
[\cite{En11}, Proposition 10.1], then $\mathrm{Hom}_R(P^0,P)$ is exact by Lemma 5.2.
Therefore, $\mathrm{Hom}_R(L^1,P)$ is exact since $\mathrm{Hom}_R(G,P)$ is so, and hence Ext$^1(L^1,P)=0$. This yields exactness of the sequence
$$0 \rightarrow \mathrm{Hom}(L^1,P) \rightarrow \mathrm{Hom}(P^0,P) \rightarrow \mathrm{Hom}(G,P) \rightarrow 0.$$
Continuously using the methods above, we get a strongly C-E exact sequence of complexes
$$0\rightarrow G\rightarrow P^0\rightarrow P^{-1}\rightarrow\cdots$$with each $P^i$ C-E projective for $i\leq0$ and which remains exact after applying Hom$(-,P)$ for any
C-E projective complex $P$.

Let $g_i:T_i\rightarrow G_i$ be a projective precover of $G_i$ for $i\in\mathbb{Z}$. Then it is clear that $T_i= G_i$ and $g_i=1_{G_i}$. We define a morphism of complexes $\varphi:\bigoplus\limits_{i\in\mathbb{Z}}D^{i}(G_i)\rightarrow G$ as follows:
$$\xymatrix{\cdots  \ar[r]^{} &  G_2\oplus G_{1} \ar[d]_{(\delta^G_{2},1)} \ar[r]^{} & G_1\oplus G_{0} \ar[d]_{(\delta^G_{1},1)}\ar[r]^{} & G_{0}\oplus G_{-1} \ar[d]_{(\delta^G_{0},1)}\ar[r]^{}& \cdots
   \\\cdots  \ar[r]^{\delta^G_{2}} & G_1 \ar[r]^{\delta^G_{1}} & G_{0}  \ar[r]^{} & G_{-1} \ar[r]^{} &\cdots.}
$$
It is clear that $\varphi$ is surjective. Put $P^1=\bigoplus\limits_{i\in\mathbb{Z}}D^{i}(G_i)$. Then we have a short exact sequence of complexes
$$0\rightarrow K^1\rightarrow P^1\stackrel{\phi}\rightarrow G\rightarrow 0\eqno(\ddag_2)$$
with $P^1$ C-E projective and $0\rightarrow (K^1)_n\rightarrow (P^1)_n\stackrel{{\phi}_n}\rightarrow G_n\rightarrow 0$ split. Obviously, $K^1$ is exact since $G$ and $P^1$ are exact, and so the sequence $(\ddag_2)$ is strongly C-E exact by Lemma 5.4. Using the methods above, we can also show that $0\rightarrow K^1\rightarrow P^1\stackrel{\phi}\rightarrow G\rightarrow 0$ is Hom$(-,P)$ exact for any C-E projective complex $P$, and so we get a strongly C-E exact sequence $$\cdots \rightarrow P^2
\rightarrow P^1
\rightarrow G \rightarrow 0 $$with $P^i$ C-E projective for $i\geq1$ and which remains exact after applying Hom$(-,P)$ for any
C-E projective complex $P$. Therefore, $G$ is strongly C-E Gorenstein projective.
\endproof

\begin{cor} Let G be an exact complex bounded right. Then G is strongly $\mathrm{C}$-$\mathrm{E}$ Gorenstein projective if and only if $G_i$ is a projective module for each $i\in\mathbb{Z}$.
\end{cor}
\proof $(\Rightarrow)$. It follows from Proposition 5.5.

$(\Leftarrow)$. It is clear by Proposition 3.8 and the fact that the class of projective modules is projectively resolving.
\endproof

As we know, an exact complex $P$ is projective if and only if $\mathrm{Z}_n(P)$ is projective in $R$-$\mathrm{Mod}$ for each $n\in\mathbb{Z}$. As a consequence of Proposition 3.8, we have the following result.

\begin{cor} Let G be an exact complex. Then G is strongly $\mathrm{C}$-$\mathrm{E}$ Gorenstein projective if and only if $G_n$ is projective and $\mathrm{Z}_n(G)$ is Gorenstein projective in $R$-$\mathrm{Mod}$ for each $n\in\mathbb{Z}$.
\end{cor}

We need the following easy lemma whose proof is routine.
\begin{lem} Let C be a complex. If $\mathrm{Hom}_R(P ,C)$ is exact for every $\mathrm{C}$-$\mathrm{E}$ projective complexes P, then C is exact.
\end{lem}

Note that projective complexes are strongly C-E Gorenstein projective complexes which are $\sharp$-projective. However, the converse is not true in general. The following examples which show that strongly C-E Gorenstein projective complexes lie strictly between projective complexes and $\sharp$-projective complexes.

\begin{exa} [[\mbox{\cite{BM07}, Example 2.5]}]Consider the quasi-Frobenius local ring $R=k[X]/(X^2)$ where $k$ is a field.  Then
$$P=:\cdots\rightarrow R\stackrel{x}\rightarrow R\stackrel{x}\rightarrow R\rightarrow\cdots$$ is a complete projective resolution in R-$\mathrm{Mod}$ and $\mathrm{Z}_i(P)$ is not projective in R-$\mathrm{Mod}$. So P is a strongly $\mathrm{C}$-$\mathrm{E}$ Gorenstein projective complex by Corollary 5.7 and P is not a projective complex.
\end{exa}

\begin{exa}  The complex A in [\cite{HJ09}, Example 2.4] is $\sharp$-projective but it is not strongly $\mathrm{C}$-$\mathrm{E}$ Gorenstein projective.
\end{exa}

We conclude this section with the following result which establishes a relationship between the above complexes.

\begin{thm} Let $G$ be a complex. Then the following statements are equivalent$:$

\noindent$(1)$  $G$ is exact, $G_n$ is projective and $\mathrm{Z}_n(G)$ is Gorenstein projective in $R$-$\mathrm{Mod}$ for each $n\in\mathbb{Z}$.

\noindent$(2)$ $G$ is strongly $\mathrm{C}$-$\mathrm{E}$ Gorenstein projective and $\mathrm{Hom}_R(P,G)$ is exact for every $\mathrm{C}$-$\mathrm{E}$ projective complex $P.$

\noindent$(3)$ $G_n$ is projective in R-$\mathrm{Mod}$, $\mathrm{Hom}_R(P,G)$ and $\mathrm{Hom}_R(G,P)$ are exact for every $\mathrm{C}$-$\mathrm{E}$ projective complex $P.$

\noindent$(4)$ $G$ is $\sharp$-projective, $\mathrm{Hom}_R(P,G)$ and $\mathrm{Hom}_R(G,P)$ are exact for every $\mathrm{C}$-$\mathrm{E}$ projective complex $P.$

\end{thm}
\proof It follows by Proposition 5.5, Corollary 5.7 and Lemma 5.8.
\endproof

Gillespie introduce and consider the notion of $\mathcal{A}$ complexes in [\cite{Gi04}].

As a direct consequence of Theorem 5.11, the following observation, which is a new characterization of $\mathcal{A}$ complexes in $\mathcal{C}(R$-Proj) when $\mathcal{A}$ is the class of Gorenstein projective modules, can be obtained. It is also a Gorenstein version of a well-known result: an exact complex $P$ with $\mathrm{Z}_n(P)$ projective for $n\in\mathbb{Z}$ if and only if $P$ is projective.

\begin{cor} Let $G$ be in $\mathcal{C}(R$-$\mathrm{Proj})$. Then $G$ is exact and $\mathrm{Z}_n(G)$ is Gorenstein projective in $R$-$\mathrm{Mod}$ for each $n\in\mathbb{Z}$ if and only if $\mathrm{Hom}_R(P,G)$ and $\mathrm{Hom}_R(G,P)$ are exact for every $\mathrm{C}$-$\mathrm{E}$ projective complex $P.$

\end{cor}

\noindent\textbf{Acknowledgements}\quad
The paper was partially written when the first named author visited Capital Normal University in 2013. The first author would like to thank Professor Changchang Xi for his warm hospitality and valuable comments on the work. The first author also would like to thank Professor Husheng Qiao for his support and our helpful discussions.


\def\refname
\end{document}